\documentclass{article}
\usepackage{amsmath,amsfonts,amssymb}

\numberwithin{equation}{section}

\newtheorem{theorem}[equation]{Theorem}
\newtheorem{lemma}[equation]{Lemma}
\newtheorem{corollary}[equation]{Corollary}
\newtheorem{proposition}[equation]{Proposition}
\newtheorem{definition}[equation]{Definition}
\newtheorem{example}[equation]{Example}

\newtheorem{remark}[equation]{Remark}
\newtheorem{remarks}[equation]{Remarks}

\newcommand{\codim}{\operatorname{codim}}

\newcommand{\Pic}{\operatorname{Pic}}

\newcommand{\Exc}{\operatorname{Exc}}

\newcommand{\PP}{{\mathbb P}}

\newcommand{\ZZ}{{\mathbb Z}}

\newcommand{\ci}{C^\prime}
\renewcommand{\o}{\omega}
\renewcommand{\O}{{\mathcal O}}
\newcommand{\I}{{\mathcal I}}
\newcommand{\F}{{\mathcal F}}
\newcommand{\E}{{\mathcal E}}
\newcommand{\G}{{\mathcal G}}
\newcommand{\Hom}{{\mathcal H}om}
\newcommand{\Mod}{{\mathcal M}od}
\newcommand{\Ext}{{\mathcal E}xt}

\newenvironment{pf}
{\noindent\textbf{Proof.}}
{\hfill{$\square$}\medskip}

{\hfill\textbf{ Q.E.D.}\medskip}

\begin{document}

\title{Linkage on singular
rational
normal surfaces and $3$-folds with application to the classification
of curves of maximal genus.}
\author{Rita Ferraro}

\maketitle
\begin{center}
Dipartimento di Matematica \\
Universit\`a di Roma  Tre\\
Largo San Leonardo Murialdo, 1 - 00146 Roma \\
\end{center}

\section{Introduction}

In algebraic geometry and commutative algebra the notion of linkage
by a complete intersection, which we will here call {\it classical linkage},
has  been for a long time an interesting and active topic.
In this note we provide a generalization of classical linkage
in  a different
context. Namely we will look at {\it 
residuals
in the scheme theoretic  intersection of a rational normal surface or
$3$-fold with two hypersurfaces  of degree $a$ and $b$}
 (a c.i. of type $(a, b)$
on the scroll, see Def. \ref{c.i}). When the scroll is singular a c.i. of
type
$(a,  b)$ on it may  not be
Gorenstein, i.e. its dualizing sheaf may   not be invertible.
If this is the case, classical linkage, even if suitably generalized, 
does not apply.

The main purpose of this article is to establish a framework allowing one 
to find relations between the dimension of some important cohomology groups 
attached to the two linked schemes.
In the last part of the paper
we  show how to apply these results and  techniques to the 
classification of curves $C$ in $\PP^n$ of degree $d$ and maximal genus 
$G(d, n, s)$ among
those not
contained in surfaces of degree less than a certain fixed one $s$.
This was the original motivation of this work.
A complete classification theorem has been given for $n=3$ by L. Gruson
and C. Peskine in 
\cite{gp}, for
$n=4$ by L. Chiantini and C. Ciliberto in \cite{cc} and for $n=5$
by the author in  \cite{f3}. For $n=3$ and $n=4$
the respective
classification Theorems have been proven with techniques of classical
linkage but for $n\geq 5$ this is no longer possible.
For $n\geq 5$ and $s\geq 2n-1$ the classification procedure consists {\it in
the precise
description of
the linked curve  to $C$ by a certain
c.i. on a
rational
normal $3$-fold $X$}. In
Example
\ref{esempio1} we describe this linked curve in the  easiest  case, i.e.
when it is a plane curve.
In Example \ref{esempio2} we construct examples of smooth curves of
maximal genus $G(d, n, s)$
for  every $d$ and $s$ in the range of Example \ref{esempio1}.

Turning to a detailed presentation of the results, 
let $W\subset\PP^{n-1}$ be a rational normal surface
and let $X\subset\PP^n$ be a rational normal $3$-fold in $\PP^n$;
throughout
the article  $W$ will be often  a general hyperplane section of $X$.
Let $Z_1$ and $Z_2$ (resp. $Y_1$ and $Y_2$) be the  two linked schemes by
a
c.i. of type $(a, b)$ on $W$ (resp. $X$).
We begin with the case when  the scroll $W$ (resp. $X$) is 
{\it smooth}. In this case  one can use a straightforward generalization
of
classical linkage
(in particular of Prop. 2.5 of \cite {ps}). Namely
the  construction
through the {\it mapping cone} of a locally free resolution 
of $\O_{Z_2}$ (resp.
$\O_{Y_2}$)
in
$\Mod(W)$ (resp. $\Mod(X)$) from a locally free resolution of $\O_{Z_1}$
(resp. $\O_{Y_1}$),
allows us to find the following results (Prop. \ref{linkss} and Prop.
\ref{linkst}):
\[
h^0(\I_{{Z_2}|W}\otimes
\O_W((i+2)H +K_W)=h^1(\I_{{Z_1}|W}\otimes\O_W((a+b-2-i)H))
\]
for $i<\min\{a, b\}$, and
\[
h^1(\I_{{Y_2}|X}\otimes\O_X((i+3)H +K_X))=
h^1(\I_{{Y_1}|X}\otimes\O_X((a+b-3-i)H))
\]
for every $i$. Here  $H$ is a hyperplane section
divisor (in both $W$ and $X$) and $K_W\sim -2H+(n-4)R$ 
(resp. $K_X\sim -3H+(n-4)R$) is the canonical 
divisor in $W$ (resp. in $X$), where $R$ is a divisor in the ruling
of $W$ (resp. $X$).

The first result allows us to compute $h^0(\I_{{Z_2}|W}\otimes
\O_W((i+2)+K_W))$ for low values of $i$ in terms of the Hilbert
function $h_{Z_1}(a+b-2-i)$ of the $0$-dimensional scheme  $Z_1$.
If $Y_1$ is aritmetically Cohen Macaulay, the second result
implies that $h^1(\I_{{Y_2}|X}\otimes\O_X((i+3)H+ K_X))=0$ for every
$i$, and therefore the restriction map
$H^0(\I_{{Y_2}|X}\otimes\O_X(iH +(n-4)R))\to
H^0(\I_{{Z_2}|W}\otimes\O_W(iH+ (n-4)R))$ is surjective for every $i$;
that means we can lift curves on $W$ linearly equivalent to
$iH+(n-4)R$ through a general hyperplane section  $Z_2$ of $Y_2$ to
surfaces on $X$ linearly
equivalent to $iH+(n-4)R$ passing through $Y_2$.
\medskip

This construction is problematic if the scroll is singular; however
a naive approach yields some interesting, even if weaker, results.
Namely,
in Section 4  we will prove for $Y_1$ arithmetically Cohen-Macaulay
and for  $i<\min\{a, b\}$ the following inequalities:
\begin{theorem}{(see Theorems \ref{linksings} and \ref{linksingt})}
\begin{eqnarray*}
h^0(\I_{{Z_2}|W}\otimes\O_W(iH+(n-4)R))&\geq &
h^1(\I_{{Z_1}|W}\otimes\O_W((a+b-2-i)H))\\ \label{s1}
h^0(\I_{{Y_2}|X}\otimes\O_X(iH+(n-4)R))&\geq &
h^1(\I_{{Z_1}|W}\otimes\O_W((a+b-2-i)H)).
\end{eqnarray*}
\end{theorem}
In this case $\O_W(iH+(n-4)R)$ (resp. $\O_X(iH+ (n-4)R)$)
is the
{\it divisorial sheaf} associated to a Weil divisor $\sim iH+(n-4)R$ on
$W$ (resp. $X$). This means that the Hilbert function
$h_{Z_1}(a+b-2-i)$ gives us a lower bound for both the dimensions
of the vector spaces of curves on $W$ passing through $Z_2$ and
surfaces on $X$ passing through $Y_2$ linearly equivalent to
$iH+(n-4)R$, for low values of $i$. As we will see, this is sufficient 
for many applications.

The same technique used to prove the above results
allows us to prove also  a formula which relates   the
arithmetic genera  of the  curves $Y_1$ and $Y_2$, linked by
a complete intersection $Y$ of type $(a, b)$ on the scroll $X$,
in the case that the vertex of $X$ is a point:
\begin{theorem}{(see Theorem \ref{genussing})}
\[
p_a(Y_2)=p_a(Y_1)-p_a(Y)+(a+b-3)\cdot\deg(Y_2)+
(n-4)\cdot\deg(R\cap Y_2)+1.
\]
\end{theorem}
(The same formula
is easily
proved with classical
linkage techniques in the smooth case, see Prop. \ref{genuss}).
\medskip

Much of these linkage techniques and the classification
for curves of maximal genus $G(d, n, s)$ in
case
$n=5$ appeared as part  of my doctoral dissertation \cite{f1}. The author
 thanks her advisor Ciro Ciliberto.

The paper has been written while the author was supported by a INDAM
scholarship.

\section{Preliminaries}

In this section, we collect  the definitions and notation to be
used in this paper, and state some of the basic results of linkage theory.
We introduce the  definitions  of geometric and algebraic
 linkage by a projective scheme $Y$, without
supposing  $Y$ to be a complete intersection.
Most of this material is well known; sometimes, however,
due to lack of a reference, we will indicate a proof.
As our primary tool is the theory of locally free resolution of sheaves,
we will include a short discussion of this, indicating
the main results we will use.
Moreover we will briefly introduce rational normal scrolls, in particular
what  we need about Weil divisors on them, including linkage.

\begin{definition}
\label{geomlink}
Let $Y_1$, $Y_2$, $Y$ be subschemes of a projective space $\PP$, then
$Y_1$
and $Y_2$ are geometrically linked by $Y$ if:
\begin{enumerate}
\item $Y_1$ and $Y_2$ are equidimensional,  have no embedded components
and  have no common components
\item $Y_1\cup Y_2=Y$, scheme theoretically.
\end{enumerate}
\end{definition}

The following Proposition is essentially Prop. 1.1 of \cite{ps}.
\begin{proposition}
\label{geomisalg}
Let $Y_1$ and $Y_2$ be closed subschemes of $\PP$ geometrically linked by
$Y$, then:
\begin{eqnarray*}
\I_{{Y_1}|Y}&\cong&\Hom_Y (\O_{Y_2}, \O_Y)\\
\I_{{Y_2}|Y}&\cong&\Hom_Y (\O_{Y_1}, \O_Y)
\end{eqnarray*}
\end{proposition}
\begin{pf}
By \cite{ps} Prop. 1.1
 we have that $\I_{{Y_1}|Y}\cong \Hom_\PP (\O_{Y_2}, \O_Y)$
and
$\I_{{Y_2}|Y}\cong \Hom_\PP (\O_{Y_1}, \O_Y)$.
Since $Y_1$ and $Y_2$ are both subschemes of $Y\subset \PP$
these isomorphisms can be rewritten as in the statement.
\end{pf}

\begin{definition}
\label{alglink}
Let $Y_1$, $Y_2$ be  projective schemes, then $Y_1$
and $Y_2$ are algebraically linked by a projective scheme $Y$ containing
them, if:
\begin{enumerate}
\item $Y_1$ and $Y_2$ are equidimensional and have no embedded
components
\item
\begin{eqnarray*}
\I_{{Y_1}|Y}&\cong &\Hom_Y (\O_{Y_2}, \O_Y)\\
\I_{{Y_2}|Y}&\cong &\Hom_Y (\O_{Y_1}, \O_Y)
\end{eqnarray*}
\end{enumerate}
\end{definition}

\begin{remark}
\label{rem2}
Let $Y\hookrightarrow\PP$ be a projective embedding of $Y$. Condition 2 of
Def. \ref{alglink} is equivalent to the following one in
every open set
$U\subset\PP$
\begin{eqnarray*}
\I_{Y}(U) : \I_{Y_2}(U)&\cong &\I_{Y_1}(U)\\
\I_{Y}(U) : \I_{Y_1}(U)&\cong &\I_{Y_2}(U)
\end{eqnarray*}
or, equivalently, we can say that
$\I_{Y_1}(U)$ (resp. $\I_{Y_2}(U)$) is the biggest ideal in $U$
such that $\I_{Y_1}(U)\cdot \I_{Y_2}(U)\subset \I_{Y}(U)$
(resp. $\I_{Y_2}(U)\cdot \I_{Y_1}(U)\subset \I_{Y}(U)$).
\end{remark}

\begin{remarks}
\label{rem3}
If $Y_1$ and $Y_2$ are geometrically linked by $Y$, then by Prop.
\ref{geomisalg} they are also algebraically linked. Moreover
if $Y_1$ and $Y_2$ are algebraically linked by $Y$ and have no common
components, then they are geometrically linked. In this case in fact we have
$I_{Y_1}(U)\cdot I_{Y_2}(U)=I_{Y_1}(U)\cap I_{Y_2}(U)$ for every open set 
$U\subset \PP$ and this implies by the previous Remark
 that
$I_{Y}(U)=I_{Y_1}(U)\cap I_{Y_2}(U)$, i.e. $Y$ is 
the scheme theoretic union of $Y_1$ and $Y_2$.
\end{remarks}

\begin{definition}
Let $Y$ be  a projective scheme. The dualizing sheaf of $Y$ is
\[
\o_Y :=\Ext^r_\PP (\O_Y, \o_\PP)
\]
where $Y\hookrightarrow \PP$ is a projective embedding of $Y$ and
$r=\codim(Y, \PP)$.
\end{definition}

For a proof of the following Theorem the reader may consult
 \cite{e} Th. 21.15 or for more details \cite{f1} Cor. 1.2.3.

\begin{theorem}
\label{ext}
Let $Y$ and $X$ be two equidimensional projective locally Cohen-Macaulay
schemes. Suppose $Y\subset X$ and let $r^\prime$ be $\codim(Y, X)$.
Let $\F$ be a sheaf in $\Mod (Y)$. Then, for every $j\geq 0$:
\[
\Ext^j_Y (\F, \o_Y)\cong \Ext^{{r^\prime}+j}_X(\F, \o_X).
\]
In particular:
\[
\o_Y\cong\Ext^{r^\prime}_X (\O_Y, \o_X).
\]
\end{theorem}

\begin{corollary}
\label{cor1}
If $Y_1$ and $Y_2$ are projective schemes algebraically linked by a
projective Gorenstein scheme $Y$,
then:
\begin{eqnarray*}
\I_{{Y_1}|Y}\otimes\o_Y &\cong &\o_{Y_2}\\
\I_{{Y_2}|Y}\otimes\o_Y&\cong &\o_{Y_1}
\end{eqnarray*}
\end{corollary}

\begin{pf}
Since $\o_Y$ is invertible we have
\[
\I_{{Y_1}|Y}\otimes\o_Y \cong\Hom_Y(\O_{Y_2}, \O_Y)\otimes\o_{Y}\cong
\Hom_Y(\O_{Y_2}, \o_Y).
\]
By Theorem \ref{ext}
$
\Hom(\O_{Y_2}, \o_Y)\cong\Ext^r_\PP(\O_{Y_2}, \o_\PP)\cong\o_{Y_2}
$
where $Y\hookrightarrow \PP$ is a projective embedding of $Y$ and
$r=\codim(Y,
\PP)$.
\end{pf}
\medskip

It is  well known (see e.g. \cite{m} \S 1.2) that, if 
$Y\subset\PP^n=\operatorname{Proj}(S)$
is a projective equidimensional locally Cohen Macaulay scheme of codimension
$r$, then there exists
a locally free resolution $\F^\bullet$ of $\O_Y$ in $\Mod(\PP^n)$ of the
type:
\begin{equation}
\label{resinP}
0\to\mathcal{K}\to \F_{r-1}\to\cdots\to \F_1\to\O_{\PP}\to\O_Y\to 0
\end{equation}
where every $\F_i$ is a finite direct sum of invertible
sheaves $\bigoplus \O_\PP(a_i)$ for 
$a_i\in \ZZ$, and
$\mathcal{K}$ is locally free.
The existence of this resolution can be proved starting from a {\it minimal 
free resolution} $F^\bullet$ for the {\it saturated ideal}
$\operatorname{{I}_{Y|\PP}}$, as a graded module over $S$, of the form:
\[
0\to
F_n\stackrel{f_n}{\rightarrow}F_{n-1}\stackrel{f_{n-1}}{\rightarrow}
\cdots\to
F_1\stackrel{f_1}{\rightarrow} S\to  {S}/{\operatorname{{I}_{Y|\PP}}}\to 0.
\]
Where the saturated ideal associated to $Y\hookrightarrow\PP^n$ is 
\[
\operatorname{{I}_{Y|\PP}}=H^0_*(\I_{Y|\PP})
:=\bigoplus_{d\in\ZZ}H^0(\PP, \I_{Y|\PP}).
\]
It then follows from the local version of the Auslander Buchsbaum Theorem 
(see \cite{e} Th. 19.9 and Cor. 19.15)
that $K_i:=\ker (f_{i-1})$ is locally free as graded $S$-module
for $i\geq r$. Therefore we can sheafify the following locally free
resolution:
\[
0\to K_{r}\to F_{r-1}\to\cdots\to F_1\to S\to 
{S}/{\operatorname{I}_{Y|\PP}}\to 0.
\]
and we obtain (\ref{resinP}).

For the purposes of this article we need a slight generalization of 
the above fact. More precisely we will use the following result: 

\begin{lemma}
\label{locfreeres}
Let $X=\operatorname{Proj} (S)$ be a projective scheme, such that $S$ 
is generated by $S_1$ as an $S_0$-module.
Let $Y\subset X$ be a projective equidimensional
 locally Cohen Macaulay scheme, contained in the smooth part of $X$, and
let $r=\codim(Y, X)$. Then:
\begin{enumerate}
\item
 There exists a locally free resolution
$\F^\bullet$ of
$\O_Y$ in $\Mod(X)$ of the type:
\[
0\to\mathcal{K}\to \F_{r-1}\to\cdots\to \F_1\to\O_X\to\O_Y\to 0
\]
where every $\F_i$ is a finite direct sum of invertible
sheaves $\bigoplus \O_X(a_i)$ for 
$a_i\in \ZZ$, and
$\mathcal{K}$ is locally free in $\Mod(X)$
\item
The above resolution can be obtained by sheafification, 
starting from a minimal free resolution $F^\bullet$ of the saturated ideal
$\operatorname{I_{Y|X}}=H^0_*(\I_{Y|X}):=\bigoplus_{d\in \ZZ}H^0(X,
\I_{Y|X})$ as a graded module over $S$.
\end{enumerate}
\end{lemma}

\begin{pf}
If $S$ is generated in degree $1$, then by \cite{h2} Prop. 5.15 there 
is a natural isomorphism ${H^0_*(\F)}^\sim \cong\F$,
 for every quasi-coherent
 sheaf $\F$
on $X$.
Therefore, one can start from a free resolution $F^\bullet$ of the
saturated ideal $\operatorname{\I_{Y|X}}$ as a graded $S$-module, then 
using exactly the same arguments as in the case $X=\PP$, one obtaines, by
sheafification,  the required resolution for $\O_Y$ in $\Mod(X)$.
\end{pf}

\begin{remark}
The existence of a locally free resolution for $\O_Y$ in $\Mod(X)$
of the type described in Lemma \ref{locfreeres} 1. can be proved also
without the assumption that $S$ is generated in degree $1$. But in this
case it cannot be obtained by sheafification from a free resolution of
the saturated ideal $\operatorname{I_{Y|X}}$. 
For a proof of the existence of
the resolution in this more general situation the reader may consult
\cite{f1}
\S 1.2.
\end{remark}

\begin{corollary}
\label{dualres}
In the hypotheses of Lemma \ref{locfreeres} and  assuming that $X$ is
locally
Cohen-Macaulay, then applying the functor $\Hom_X(\cdot, \o_X)$
to $\F^\bullet$ we  obtain a resolution of lenght $r$ of $\o_Y$ in
$\Mod(X)$.
\end{corollary}
\begin{pf}
It follows from the fact that $\Ext^i(\F, \o_X)=0$ for every $i<r$ and
every coherent sheaf $\F$ in $\Mod(Y)$, and from Theorem \ref{ext}. For
details see
\cite{f1} Prop. 1.2.4. To prove that $\Ext^i(\F, \o_X)=0$ for every $i<r$
one can  prove, as in \cite{h2} III Lemma 7.3,  that $\Ext^i(\F,
\o_\PP)=0$ for every
$i<\codim(Y,
\PP)$, where $X\hookrightarrow \PP$, then use Th. \ref{ext}. 
\end{pf}

The next  Proposition  is a little generalization of
Prop.
2.5 of \cite{ps}, therefore we omit the proof. 
\begin{proposition}
\label{mapcon} 
Let $Y_1$ and $Y_2$ be two equidimensional projective schemes
algebraically linked by a projective Gorenstein scheme  $Y$.
Suppose that $Y_1$ is locally Cohen Macaulay and suppose that $Y$
is contained in the smooth part of a locally Cohen Macaulay
projective scheme $X=\operatorname{Proj}(S)$; let $r=\codim(Y, X)$. Let
$\F^\bullet$ and
$\F^\bullet_1$ be locally free resolutions in $\Mod (X)$ of $\O_Y$
and $O_{Y_1}$ respectively, as in Lemma \ref{locfreeres}. Let
$r^.:\F^\bullet \to \F^\bullet_1$ be  a morphism of complexes
induced by the restriction map $r: \O_Y \to \O_{Y_1}$, let the
dual functor $\vee$ be the functor $\Hom_X(\cdot, \O_X)$, then the
mapping cone
$C^\bullet ({r^.}^\vee)$ of the morphism ${r^.}^\vee :
{\F_1^\bullet}^\vee \to {\F^\bullet}^\vee$ is a locally free
resolution of lenght $r$ of ${\o_{Y}}_{|{Y_2}}\otimes {\o_X}^\vee$
in $\Mod(X)$.
\end{proposition}

However we stress that  the hypothesis that $S$ is generated in degree 
$1$ is essential for the existence of the morphism $r^.:\F^\bullet \to
\F^\bullet_1$. Indeed, in this case, $r^.$ is induced, by sheafification,
from the morphism ${\rm r}^.:F^\bullet\to F^\bullet_1$ of the free
resolutions of $S/\operatorname{I_{Y}}$ and $S/\operatorname{I_{Y_1}}$,
as graded modules over $S$.

\medskip

\begin{definition}
\label{c.i}
Let $X$ be a projective scheme of dimension $k$; let $a_i\in \mathbb{N}^+$
and
let $1\leq r\leq k$. A complete intersection (shortly c.i.)
on $X$ of kind $(a_1, \dots, a_r)$ is an equidimensional projective scheme
$Y\subset
X$ such that $\codim(Y, X)=r$, which is scheme theoretic intersection of
Cartier divisors $D_i\in |\O_X(a_i)|$ for $i=1, \dots, r$.
\end{definition}

\begin{corollary}
\label{c.i.res}
In the hypotheses of Prop. \ref{mapcon} if $Y\subset X$ is a c.i.
 on $X$ of kind $(a_1, \dots, a_r)$, then $C^\bullet
{({r^.}^\vee)}(-a_1\cdots -a_r)$ is a locally free resolution of
$\O_{Y_2}$ in $\Mod(X)$.
\end{corollary}
\begin{pf}
$Y$ is contained in the smooth part of $X$, so we can
think of $\o_X$ as an invertible sheaf on $Y$; by adjunction formula
$\o_Y$ is $\O_Y\otimes\o_X(a_1+\cdots +a_r)$.
\end{pf}

\medskip

We want now to fix some notation about rational normal scrolls
and point out
what we will need in the next sections.
A rational normal scroll $X\subset \PP$ of dimension $r$ and degree
$f$  is the
image of a projective bundle $\PP(\E)\to \PP^1$ over $\PP^1$
through the morphism $j$ defined by the tautological line bundle
$\O_{\PP(\E)}(1)$, where $\E\cong \O_{\PP^1}(a_1)\oplus \cdots \oplus
\O_{\PP^1}(a_{r})$ with $0\leq a_1\leq \cdots \leq a_{r}$ and
$\sum a_i=f=n-r$. If $a_1=\cdots=a_l=0$, $1\leq l <r$, $X$ is
singular and the  vertex $V$ of $X$ has dimension $l-1$. Let us
denote $\PP(\E)=\tilde X$. The morphism  $j: \tilde{X}\to X$ is a
rational resolution of singularities, i.e. $X$ is {\it normal} and
{\it arithmetically Cohen-Macaulay} and $R^ij_*\O_{\tilde X}=0$
for $j>0$. We will call $j:\tilde X\to X$ the {\it canonical
resolution} of $X$.  It is well known that
$\Pic(\tilde{X})=\mathbb{Z}[\tilde{H}]\oplus\mathbb{Z}[\tilde{R}]$,
where $[\tilde{H}]=[\O_{\tilde{X}}(1)]$ is the hyperplane class
and $[\tilde{R}]=[\pi^*\O_{\PP^1}(1)]$ is the class of the fibre
of the  map $\pi: \tilde X\to \PP^1$. The intersection form on $\tilde X$
is determined by the rule: 
\[
{\tilde H}^{r}=f \qquad
{\tilde H}^{r-1}\cdot {\tilde R}=1 \qquad
{\tilde H}^{r-2}\cdot {\tilde R}^2=0 
\]
Let us denote with $X_S$
the smooth part of $X$ and with $\Exc(j)$ the exceptional locus of $j$.
Then
$j:\tilde X\setminus \Exc(j) \to X_S$ is an isomorphism. 
Let
$H$ and $R$ be the {\it strict images} of $\tilde H$  and $\tilde
R$
respectively (i.e. the {\it scheme theoretic closure} 
$\overline {j({\tilde H}_{|j^{-1}X_S})}$
and
$\overline {j({\tilde R}_{|j^{-1}X_S})}$). Then   we have the following
well known result:
\begin{lemma}
\label{Cl} 
Let $X\subset \PP^n$ be a  rational normal scroll of
degree $f$  and let $j:\tilde X\to X$ be its canonical resolution.
Let $\operatorname{Cl}(X)$ be the group of Weil divisors on $X$
modulo linear equivalence. Then
\begin{enumerate}
\item
If $\codim (V, X)>2$, $\operatorname{Cl}(X)\cong\ZZ [H]\oplus\ZZ
[R]$;
\item
If $\codim(V, X)=2$, $H\sim f R$ and
$\operatorname{Cl}(X)\cong\mathbb{Z}[R]$.
\end{enumerate}
\end{lemma}

We recall here from \cite{f2} the definition of {\it proper} and 
{\it (integral)
total transform}
of a Weil divisor in $X$.
In the last section (Example \ref{esempio1}) we will use  proper and
integral total transforms together with \cite{f2} Prop.
4.11 to compute
the multiplicity of the vertex $V$ in the  intersection scheme of two
effective divisors on a
rational normal $3$-fold  $X$ with $\codim(V, X)=2$. 

\begin{definition}
\label{proper}
Given a prime divisor $D$ on $X$, the  proper transform $\tilde D$ of
$D$ in $\tilde X$ is the scheme theoretic closure $\overline{j^{-1}(D\cap
X_S)}$. The proper transform of any Weil
divisor in
$X$ is then defined by linearity.
\end{definition}
\begin{definition}
\label{totale}
Let $\codim(V, X)=2$ and let $D\sim dR$ be an effective Weil divisor
on $X$, divide $d-1=kf+h$ ($k\geq -1$ and $0\leq h < f$), then the
{integral total transform} of $D$ in $\tilde X$ is $D^*\sim (k+1)\tilde H
-(f-h-1)\tilde R$.
\end{definition}

Let us define  on $X$ the following coherent sheaves for $a,
b\in \mathbb{Z}$:
\begin{definition}
\[
\O_X(a, b):= j_*\O_{\tilde{X}}(a\tilde{H}+b\tilde{R}).
\]
\end{definition}
We will usually write $\O_X(a)$ instead of $\O_X(a, 0)$. Moreover
for every coherent sheaf $\F$ on $X$ we will write $\F(a, b)$ instead
of $\F\otimes\O_X(a, b)$.
If the scroll $X$ is smooth, then the sheaves $\O_X(a, b)$ are
the invertible sheaves associated to the Cartier divisors $\sim aH+bR$ 
while when $X$ is singular this is no longer true.
In this case we have the following Proposition which is proved in
\cite{f2}
 (Cor. 3.10 and Th. 3.17). The reader may refer to \cite{h1} for a survey
on {\it divisorial sheaves} associated to generalized divisors.

\begin{proposition}
\label{divsheaf}
Let $X\subset \PP^n$ be a singular rational normal scroll of degree $f$,
dimension $r$
and vertex $V$, then:
\begin{enumerate}
\item
If $\codim(V, X)>2$ the sheaf $\O_X(a, b)$ is reflexive for every 
$a, b\in\ZZ$ and it is the divisorial sheaf associated to  a Weil divisor
$\sim aH+bR$;
\item
If $\codim(V, X)=2$ the sheaf $\O_X(a, b)$ is reflexive for every 
$a, b\in\ZZ$ such that $b<f$; in this case the sheaves 
$\O_X(a, b)$ with $a+fb=d$ are all isomorphic to the 
the divisorial sheaf associated to a  Weil divisor
$\sim dR$;
\end{enumerate} 
\end{proposition}

In the hypoteses  of 
Prop. \ref{divsheaf}, the dualizing
sheaf $\o_X$ of $X$ is (see \cite{s}):
\begin{equation}
\label{dual}
\o_X=j_*\O_{\tilde X}(K_{\tilde X})=\O_X(-r, f-2).
\end{equation}
By Prop. \ref{divsheaf}, $\o_X$ is a divisorial sheaf, 
therefore the {\it canonical divisor} of $X$ is
$K_X\sim -rH+(f-2)R$.


The following result is essentially due to Hartshorne
(Linkage of generalized divisors by a complete intersection:
\cite{h1}, Prop. 4.1). He states it for divisors on a complete
intersection but the same proof goes over as well.

\begin{proposition}{(Linkage of divisors)}
\label{linkdiv}
Let $D_1$ be an effective Weil divisor on a rational normal scroll
$X\subset \PP^n$. Let $F\subset\PP^n$ be a hypersurface containing $D_1$;
let
$D$
be the Cartier divisor on $X$ defined by $F$, then the effective divisor
$D_2=D-D_1$ is algebraically linked to $D_1$ by $D$.
\end{proposition}
\begin{pf}
By Definition \ref{alglink} of algebraic linkage we should prove the
isomorphisms
$
\I_{{D_1}|D}\cong\Hom_D(\O_{D_2}, \O_D)
$ and
$
\I_{{D_2}|D}\cong\Hom_D(\O_{D_1}, \O_D)$;
 by simmetry it is enough to prove
the second formula.
First we apply  the functor $\Hom_X(\O_{D_1}, \cdot)$ to the exact
sequence
\[
0\to \I_{D|X}\stackrel{\alpha}\to\O_X\to \O_D\to 0.
\]
This gives
\[
0\to \Hom_X(\O_{D_1}, \O_D)\to 
\Ext^1_X(\O_{D_1}, \I_{D|X})\stackrel{{\alpha^\prime}}
\to\Ext^1_X(\O_{D_1}, \O_X).
\] 
Since $D$ is a Cartier divisor on $X$, the map $\alpha$ is locally
 multiplication by a non-zero divisor $f$. Since $D_1\subset D$, this
element $f$ annihilates $\O_{D_1}$, and so the induced map
$\alpha^\prime$
on $\Ext$ is zero. Hence we find
\[
\Hom_X(\O_{D_1}, \O_D)\cong\Ext^1_X(\O_{D_1}, \I_{D|X}).
\]
Next we apply  the functor $\Hom_X(\cdot, \I_{D|X})$ to the exact
sequence
\[
0\to \I_{{D_1}|X}\to\O_X\to \O_{D_1}\to 0.
\]
This gives
\[
0\to \Hom_X(\O_{X}, \I_{D|X})\to 
\Hom_X(\I_{{D_1}|X}, \I_{D|X})\to\Ext^1_X(\O_{D_1}, \I_{D|X})\to 0.
\] 
Here the first term is just $\I_{D|X}$. Since $\I_{D|X}$ is invertible,
the
second term is $\I_{D|X}\cdot{\I_{{D_1}|X}}^\vee=\I_{{D_2}|X}$ because
$D_2=D-D_1$.
The third term we have identified above. Since the quotient of
$\I_{D_2|X}$
by $\I_{D|X}$ is $\I_{D_2|D}$ we obtain the desired formula
\[
\I_{{D_2}|D}\cong\Hom_X(\O_{D_1}, \O_D).
\]
\end{pf}

\section{Linkage by a complete intersection in the smooth case}

In this section $W$ is a smooth rational normal surface 
in
$\PP^{n-1}$ and
$X$ is a smooth rational normal $3$-fold in $\PP^{n}$, according to
(\ref{dual})
the dualizing
sheaves of $W$ and $X$ are respectively
$\omega_W=\O_W(-2, n-4)$ and $\omega_X=\O_X(-3, n-4)$.
Since we are in the smooth case, we consider algebraic linkage, where the
subschemes
need not have distinct components.

\begin{proposition}
\label{linkss}
Let $W\subset\PP^{n-1}$ be a smooth rational normal surface.
Let $Z_1\subset W$ be  a projective $0$-dimensional
locally Cohen-Macaulay scheme. Let $Z_2\subset W$ be a
projective $0$-dimensional scheme. Let $Z=W\cap F_a\cap F_b$ be a
c.i. on $W$ of type $(a, b)$.
Assume that $Z_1$ and $Z_2$ are algebraically linked
by $Z$.
Then for $i<\min\{a, b\}$:
\begin{equation}
\label{eq1linkss}
h^0(\I_{{Z_2}|W}(i, n-4))=h^1(\I_{{Z_1}|W}(a+b-2-i))
\end{equation}
or equivalently, in terms of the Hilbert function of $Z_1$:
\begin{equation}
\label{eq2linkss}
h^0(\I_{{Z_2}|W}(i, n-4))=\deg Z_1-h_{Z_1}(a+b-2-i).
\end{equation}
\end{proposition}
\begin{pf}
Since $W$ is smooth and $Z_1$ is locally Cohen-Macaulay with
$\codim (Z_1, W)=2$, by Lemma \ref{locfreeres} 
there exists a locally free resolution of $\O_{Z_1}$
in $\Mod(W)$ that  looks like:
\begin{equation}
\label{z1}
0\to \mathcal{F}\to F\to \O_W \to \O_{Z_1}\to 0
\end{equation}
where $F$ is a finite direct sum $\bigoplus\O_W(\alpha_i)$ with
$\alpha_i\in\ZZ$ and
$\mathcal{F}$ is locally free.

The resolution of $\O_Z$ in $\Mod(W)$ is:
\begin{equation}
\label{z}
0\to \O_W (-a-b)\to \O_W (-a)\oplus\O_W (-b)\to \O_W \to \O_{Z}\to 0
\end{equation}

Let $r: \O_Z\to \O_{Z_1}$ be the morphism of sheaves induced by
$Z_1\subset Z$.
By Cor. \ref{c.i.res} the mapping cone  $C^\bullet({r^.}^\vee)(-a-b)$
is a locally free resolution of $\O_{Z_2}$ in
$\Mod(W)$:
\begin{equation}
\label{z2}
0\to F^\vee (-a-b)\to \O_W (-a)\oplus\O_W (-b)\oplus\mathcal{F}^\vee (-a-b)
\to \O_W \to \O_{Z_2}\to 0
\end{equation}
We look now at the new exact sequence obtained by tensoring \ref{z2}
 with the invertible sheaf $\O_W (i, n-4)$. By Serre's
duality we have $h^1(F^\vee (i-a-b, n-4))=h^1(F(a+b-2-i))=0$,
since $W$ is aritmetically Cohen-Macaulay and $F=\bigoplus\O_W
(\alpha_i)$. Since $h^0(\O_W(\alpha, \beta )=0$ for $\alpha<0$,
we get for $i<\min\{a, b\}$:
\begin{equation}
\label{h0z2}
h^0(\I_{{Z_2}|W}(i, n-4)) =h^0({\mathcal F}^\vee(i-a-b, n-4))
-h^0(F^\vee(i-a-b, n-4)). 
\end{equation}

From the exact sequence
\[
0\to\I_{{Z_1}|W}\to\O_W\to\O_{Z_1}\to 0
\]
tensored with $\O_W(a+b-2-i)$ we get
$h^2(\I_{{Z_1}|W}(a+b-2-i)=h^2(\O_W(a+b-2-i))=
h^0(\O_W (i-a-b, n-4))=0$, for $i<a+b$.
After tensoring (\ref{z1})  by
$\O_W(a+b-2-i)$, we obtain for $i<a+b$:
\begin{align}
h^1(\I_{{Z_1}|W}(a+b-2-i))&=h^2({\mathcal F}(a+b-2-i))
-h^2(F(a+b-2-i)) \nonumber \\
&= h^0(\mathcal{F}^\vee(i-a-b, n-4))\label{h1z1}\\
&-h^0(F^\vee(i-a-b, n-4)).\nonumber
\end{align}
Formula
(\ref{eq1linkss}) follows from (\ref{h0z2}) and (\ref{h1z1}).
Moreover, since $W$ is aritmetically Cohen-Macaulay, from the exact sequence
\[
0\to\I_{W|\PP}\to\I_{{Z_1}|\PP}\to\I_{{Z_1}|W}\to 0
\]
we get
$h^1(\I_{{Z_1}|W}(k))=h^1(\I_{{Z_1}|\PP}(k))$ for every $k$. Moreover
$h^1(\I_{{Z_1}|\PP}(k))=h^0(\O_{Z_1}(k))-h_{Z_1}(k)=
\deg(Z_1)-h_{Z_1}(k)$. This proves (\ref{eq2linkss}).
\end{pf}

In the next proposition we consider the case of a c.i. of type $(a, b)$
on a rational normal $3$-fold $X$. The proof is similar to the proof of
Prop. \ref{linkss} and therefore we omit it.

\begin{proposition}
\label{linkst}
Let $X\subset\PP^n$ be a smooth rational normal $3$-fold.
Let $Y_1, Y_2\subset X$ be  projective equidimensional
of dimension $1$ schemes.
Assume that  $Y_1$ is locally Cohen-Macaulay.
 Let $Y=X\cap
F_a\cap F_b$ be a
c.i. of type $(a, b)$ on $X$.
Assume that $Y_1$ and $Y_2$ are algebraically linked
by $Y$.
Then for every $i$:
\[
h^1(\I_{{Y_2}|X}(i, n-4))=h^1(\I_{{Y_1}|X}(a+b-3-i))
\]
\end{proposition}

From Prop. \ref{linkst} it follows easily that if we suppose $Y_1$
arithmetically Cohen-Macaulay, then we can lift
divisors $\sim iH+(n-4)R$ on a general hyperplane
section $W\subset\PP^{n-1}$
of $X$ containing the general hyperplane section $Z_2$ of $Y_2$ to
divisors $\sim iH+(n-4)R$
on $X$ containing $Y_2$. Namely we have the following Corollary.

\begin{corollary}
\label{acm}
In the hypotheses  of Proposition \ref{linkst}, if we suppose $Y_1$
to be arithmetically Cohen-Macaulay, then the map
\[
H^0(\I_{{Y_2}|X}(i, n-4))\to H^0(\I_{{Z_2}|W}(i, n-4))
\]
is surjective for every $i$.
\end{corollary}
\begin{pf}
Since both $Y_1$ and $X$ are arithmetically Cohen-Macaulay 
we have that
$h^1(\I_{{Y_1}|X}(k))=h^1(\I_{{Y_1}|\PP}(k))=0$
 for every $k$.
By Prop. \ref{linkst} we then have
$h^1(\I_{{Y_2}|X}(i, n-4))=0$ for every $i$; the statement
follows now from the exact sequence
\[
0\to\I_{{Y_2}|X}(i-1, n-4)\to\I_{{Y_2}|X}(i, n-4)\to
\I_{{Z_2}|W}(i, n-4)\to 0.
\]
\end{pf}

The next result is a formula which relates the arithmetic genera
of the  curves 
$Y_1$,
$Y_2$ and
$Y$.

\begin{proposition}
\label{genuss}
In the hypoteses of Prop. \ref{linkst} we have the following
formula, relating the arithmetic genera of the linked curves:
\begin{equation}
\label{genusform1}
p_a(Y_2)=p_a(Y_1)-p_a(Y)+(a+b-3)\cdot\deg(Y_2)+
(n-4)\cdot\deg(R\cap Y_2)+1.
\end{equation}
\end{proposition}
\begin{pf}
First we tensor by the invertible sheaf
 $\o_Y\cong\O_Y(a+b-3, n-4)$  the exact sequence
\[
0\to\I_{{Y_2}|Y}\to \O_Y\to\O_{Y_2}\to 0.
\]
By Cor. \ref{cor1}, this gives the exact sequence
\begin{equation}
\label{seq1}
0\to\o_{Y_1}\to \o_Y\to{\o_Y}_{|Y_2}\to 0.
\end{equation}
Without loss of generality we can suppose $a, b\geq 2$. The sheaf
$\O_X(a+b-3, n-4)$ is very ample and intersects $Y_2$  in
a
divisor
$D$ such that $\O_{Y_2}(D)={\o_Y}_{|Y_2}$ is an invertible sheaf.
By Riemann-Roch we then obtain
\begin{equation}
\label{rel1}
\chi(\O_{Y_2}(D))=\chi({\o_Y}_{|Y_2})=1-p_a(Y_2)+\deg(D)
\end{equation}
(where $\deg(D)$ is the lenght of the $0$-dimensional scheme $D$).
Formula (\ref{genusform1}) follows now by  
(\ref{rel1}) and (\ref{seq1}) since
$\deg(D)=(a+b-3)\cdot\deg(Y_2)+(n-4)\cdot\deg(Y_2\cap R)$.
\end{pf}

\section{Linkage by a complete intersection on  singular rational normal
surfaces and $3$-folds}

Throughout this section $W\subset\PP^{n-1}$
will be a singular rational normal surface and $X\subset\PP^n$ a
singular rational normal $3$-fold. The vertex $V$ of $W$
is a point, while the vertex $V$ of $X$ can be either a point or a line.
\begin{lemma}
\label{prel}
Let $W\subset\PP^{n-1}$ (resp. $X\subset\PP^n$) 
be a singular rational normal surface (resp. $3$-fold).
Let $A=Q_1\cap \dots \cap Q_{n-3}$ be a generic complete intersection 
of $n-3$ quadrics of $\PP^{n-1}$ (resp. of $\PP^n$)
containing $W$ (resp. in  $X$).
Let $B=A-W$ (as Weil divisor in $Q_1\cap \dots \cap Q_{n-4}$) be the
residual scheme to $W$ (resp. to $X$) by $A$. 
Let $Y_B=W\cap B$ (resp
$Y_B=X\cap B$). Then
the scheme $Y_B$ is a divisor in $W$ (resp. $X$) linearly equivalent to
$(n-4)H-(n-4)R$.
\end{lemma}
\begin{pf}
Hartshorne's Connectedness
Theorem
(\cite{e} Th. 18.12) implies that $Y_B$ has pure
codimension $1$ in $W$ (resp. in $X$), therefore it
is a divisor.
To show that $Y_B\sim (n-4)H-(n-4)R$
 we  first want to prove 
that $Y_B$ has  degree $(n-4)(n-3)$.
Taking   a general hyperplane section of $W$
(resp. a general $(n-2)$-plane  section  of $X$), we obtain a rational
normal curve
$C_{n-2}\subset\PP^{n-2}$. Let  $A_H$, $B_H$ and $Y_{B_H}$ be
general
sections in $\PP^{n-2}$ of $A$, $B$ and $Y_B$ respectively.
By Prop. 4.1 in \cite{h1},  $C_{n-2}$
and $B_H$ are algebraically linked by $A_H$. 
From the exact sequence:
\[
0\to\I_{{B_H}|{A_H}}\to \O_{A_H}\to \O_{B_H}\to 0
\]
tensored with the invertible sheaf $\omega_{A_H}\cong\O_{A_H}(n-5)$
we obtain  by Cor. 2.4 the exact sequence:
\[
0\to\omega_{C_{n-2}}\to \omega_{A_H}\to \omega_{{A_H}|{B_H}}\to 0
\]
and we compute
\begin{equation}
\label{1}
\chi (\omega_{{A_H}|{B_H}})=\chi(\omega_{A_H})-\chi(\omega_{C_{n-2}})=
p_a(A_H)-1+1=p_a(A_H).
\end{equation}
Moreover
\begin{align}
\chi
(\omega_{{A_H}|{B_H}})&=\chi(\O_{B_H}(n-5))=\chi(\O_{B_H})+(n-5)
\deg({B_H})\nonumber\\
&=1-p_a(B_H)+(n-5)(2^{n-3}-(n-2)). \label{2}
\end{align}
Since  $p_a(A_H)=1+2^{n-4}(n-5)$, from (\ref{1}) and (\ref{2}) we obtain:
\begin{equation}
\label{3}
p_a(B_H)=2^{n-4}(n-5)-(n-2)(n-5).
\end{equation}
Substituting (\ref{3}) in the Noether's formula:
\[
p_a(A_H)=p_a(B_H)+p_a(C_{n-2})+\deg (B_H\cap C_{n-2})-1
\]
 we finally get
\[
\deg (B_H\cap C_{n-2})=(n-4)(n-3).
\]
Now we want to prove that
the minimal degree of a hypersurface in $\PP^{n-1}$ ( resp.
in $\PP^n$) containing $Y_B$ but not $W$ ( resp. not $X$) is
$n-4$. In fact this forces $Y_B$  to  be linearly equivalent
to  $(n-4)H-(n-4)R$.
First we prove that an hypersurface of  degree $d<n-4$ containing $Y_B$
is forced to  contain  $W$ (resp. $X$) and this follows from
$\deg(Y_B)=(n-4)(n-3)>d\cdot\deg(W)=d\cdot(n-2)$. Then we prove that
$h^0(\I_{{Y_B}|W}(n-4))>0$ (resp. $h^0(\I_{{Y_B}|X}(n-4))>0$).
Since $W$ (resp. $X$) is arithmetically Cohen-Macaulay, from the exact
sequence 
$0\to \I_{Y_B|W}\to \O_W\to \O_{Y_B}\to 0$ 
we get
$h^1(\I_{{Y_B}|W}(k))=0$ (resp. $h^1(\I_{{Y_B}|X}(k))=0$) for every $k$.
Therefore from the exact sequence
$0\to \I_{Y_B|W}(n-5)\to \I_{Y_B|W}(n-4)\to \I_{{Y_{B_H}}|C_{n-2}}(n-4)\to
0$ 
we
get
$h^0(\I_{{Y_B}|W}(n-4))=h^0(\I_{{Y_{B_H}}|C_{n-2}}(n-4))$
(resp. $h^0(\I_{{Y_B}|X}(n-4))=h^0(\I_{{Y_{B_H}}|C_{n-2}}(n-4))$).
Moreover
$h^0(\I_{{Y_{B_H}}|C_{n-2}}(n-4))=h^0(\O_{\PP^1}(n-4))>0$.
\end{pf}
\begin{remark}
In case of $Y_B\subset W$, or $Y_B\subset X$  and $X$ is singular along a
line, we have
$Y_B\sim (n-4)(n-3)R$ 
 since in both cases $H\sim (n-2)R$ by Lemma \ref{Cl}.
\end{remark}

\begin{remark}
\label{genA}
Let $Z=W\cap F_a\cap F_b$ (resp. $Y=X\cap F_a\cap F_b$)
be a c.i.  of type
$(a, b)$ on $W$ (resp. on $X$).
By generality we can choose $A$ as in Lemma \ref{prel}
in such a way that neither $F_a$ nor $F_b$
contain any  component of $B$ and  $Y_{{B}}$
does not contain any component of $Z$ (resp. of  $Y$), with the exception
of
the
vertex
$V$ of $W$ (resp. of $X$, if the vertex of $X$ is a line)
in the case that both $F_a$ and
$F_b$ pass through it.
\end{remark}

\begin{lemma}
\label{linkAs}
Let $W\subset\PP^{n-1}$ be a singular rational normal surface and let $A$
be as in  Lemma \ref{prel}.
Let $Z_1$ be a $0$-dimensional locally Cohen-Macaulay projective
scheme in $\PP^{n-1}$.
Let $Z_3\subset \PP^{n-1}$ be the $0$-dimensional scheme algebraically
linked
to $Z_1$ by a complete intersection $Z_A$ of type $(2, \dots, 2, a, b)$.
Suppose $a, b> n-4$.
 Then,
for $i+n-4<\min\{ a, b\}$:
\begin{equation}
\label{eq1linkAs}
h^0(\I_{{Z_3}|A}(i+n-4))=h^1(\I_{{Z_1|W}}(a+b-2-i)).
\end{equation}
\end{lemma}
\begin{pf}
Let $r: \O_{Z_A}\to \O_{Z_1}$ be the restriction map induced by the inclusion
$Z_1\subset Z_A$. Starting from the locally free resolutions of lenght
$n-1$
 of $\O_{Z_1}$
and $\O_{Z_A}$ in $\Mod(\PP^{n-1})$,
then the mapping cone   $C^\bullet ({r^.}^\vee)(-2\cdots -2-a-b)$
gives a locally free resolution of lenght $n-1$ of $\O_{Z_3}$
(Prop. 2.5 \cite{ps}).
Chasing through the
locally free resolutions of $\O_{Z_A}$, $\O_{Z_3}$ and $\O_{Z_1}$
we find:
\begin{equation*}
h^1(\I_{{Z_1}|{\PP^{n-1}}}(a+b-2-i))=h^0(\I_{{Z_3}|{\PP^{n-1}}}(i+n-4))
-h^0(\I_{{Z_A}|{\PP^{n-1}}}(i+n-4)).
\end{equation*}
From the projective normality of $W\subset \PP^{n-1}$ and 
$A\subset \PP^{n-1}$  we have that
$h^1(\I_{{Z_1}|W}(k))=h^1(\I_{{Z_1}|\PP^{n-1}}(k))$ for every $k$, and
$h^0(\I_{{Y}|A}(k))=h^0(\I_{{Y}|{\PP^{n-1}}}(k))
-h^0(\I_{A|{\PP^{n-1}}}(k))$
for every projective scheme $Y\subset A$ and every $k$.
Then the previous equality can be rewritten as
\begin{equation}
\label{eq2linkAs}
h^1(\I_{{Z_1}|{W}}(a+b-2-i))=h^0(\I_{{Z_3}|{A}}(i+n-4))
-h^0(\I_{{Z_A}|{A}}(i+n-4)),
\end{equation}
 which proves the statement for $i+n-4< \min\{a, b\}$.
\end{pf}

\begin{theorem}
\label{linksings}
Let $W\subset\PP^{n-1}$ be a singular rational normal
surface.
Let $Z_1, Z_2\subset W$ be  projective $0$-dimensional schemes. Assume
that $Z_1$ is
locally Cohen-Macaulay. Let $Z=W\cap F_a\cap F_b$  be a
c.i. of type $(a, b)$ on $W$.
Assume that $Z_1$ and $Z_2$ are geometrically linked
by $Z$.
Then for $i<\min\{a, b\}$:
\begin{equation}
\label{eq1linksings}
h^0(\I_{{Z_2}|W}(i, n-4))\geq h^1(\I_{{Z_1}|W}(a+b-2-i)) 
\end{equation}
or equivalently, in terms of the Hilbert function of $Z_1$:
\[
h^0(\I_{{Z_2}|W}(i, n-4))\geq\deg Z_1-h_{Z_1}(a+b-2-i).
\]
\end{theorem}

\begin{pf}
Let $A\subset\PP^{n-1}$ be a
complete intersection of quadrics containing $W$ and let $B$ be  the
(geometrically)
linked scheme to $W$ by $A$ as in Lemma \ref{prel}.
Let $Z_A$ and $Z_3$ be as in Lemma \ref{linkAs}.
 Let $Z^\prime=Z_A-Z$ as Weil divisor 
on $A\cap F_a$, then by \cite{h1} Prop. 4.1 we know that $Z^\prime$ is the
algebraically
linked scheme to $Z$ by $Z_A$. It follows that 
$Z^\prime=B\cap F_a\cap F_b$, i.e. $Z^\prime$ is a
 complete
intersection of type $(a, b)$ on  $B$ and that $Z^\prime =Z_3-Z_2$.

The hypersurfaces of degree $i+n-4$ containing $Z_3$  contain then
$Z^\prime$; if we suppose $i+n-4< \min\{a, b\}$
these hypersurfaces are forced to contain $B$, since $Z^\prime$ is a c.i. 
of type $(a, b)$ on it. This means that every hypersurface $T$
in the linear system $|\I_{{Z_3}|A}(i+n-4)|$ cuts on $W$ a divisor which 
split
in the union of  $Y_{B}\sim (n-4)(n-3)R$
and a divisor
$D_T\sim (i+n-4)H-(n-4)(n-3)R\sim iH+(n-4)R$.

We claim that $Z_2$ is contained in the residual divisor $D_T
$.

If we choose  $A$ as in Remark \ref{genA} it is enough to prove
it locally in a open affine subset  $U$ of $\PP^{n-1}$
containing  the vertex $V$ of $W$. So, let us suppose
$Z_2$ does contain  $V$ (otherwise there is nothing to prove). Since
$Z_1$ and $Z_2$ are geometrically linked they are disjoint, therefore the restriction
of $Z_3$ on to $U$ is a complete intersection in $U$.
Let $f=0$, $g=0$ be local equations in $U$
of $F_a$ and $F_b$ respectively, and let $q_1, \dots, q_{{n-2}\choose 2}$
be local equations of $Q_1, \dots, Q_{{n-2}\choose 2}$ , where
$Q_1, \dots, Q_{{n-2}\choose 2}$ are generators in the homogeneous
ideal $I_{W|{\PP}}$ and
$Q_1, \dots, Q_{n-3}$ are generators in the homogeneous
ideal $I_{A|{\PP}}$.
Then the ideal of ${Z_2}_{|U}$ in $U$ is $\I_{{Z_2}}(U)=
(q_1, \dots, q_{{n-2}\choose 2},
f, g)$ and the one of ${Z_3}_{|U}$ is $I_{{Z_3}}(U)=(q_1, \dots, q_{n-3}, f, g)$.

Let $\I_B(U)$ and $\I_{Y_B}(U)$ be the ideals in $U$ of $B_{|U}$
and ${Y_B}_{|U}$ respectively, then the  ideal of $Z^\prime_{|U}$ in $U$
is $\I_{Z^\prime}(U)=(\I_B(U), f, g)$. 
Since $Z_2$ and $Z^\prime$ are algebraically
linked by $Z_3$, then $\I_{{Z_2}}(U)$ is the biggest ideal in $U$ (see
Remark \ref{rem2}) such that
\begin{equation}
\label{eq2linksings}
\I_{{Z_2}}(U)\cdot (\I_B(U), f, g)\subset (q_1, \dots ,q_{n-3}, f, g).
\end{equation}

Let $T$ be a hypersurface of degree $i+n-4$ in $\PP^{n-1}$ containing
$Z_3$ but not containing $W$, let $(T)$ be the divisor cut by $T$
on $W$, then the  divisor
$D_T={(T)}-Y_B$ and $Y_B$ are algebraically linked by
${(T)}$ (Prop. \ref{linkdiv}).
By Remark \ref{rem2} $\I_{D_T}(U)$ is the biggest ideal such
that
\begin{equation}
\label{eq3linksings}
\I_{D_T}(U)\cdot  \I_{{Y_B}}(U)
\subset (q_1, \dots , q_{{n-2}\choose 2}, t),
\end{equation}
where $t$ is the local equation of $T$ in $U$.
We want now to prove that on $U$
\begin{equation}
\label{eq4linksings}
\I_{D_T}(U)\cdot  (\I_{B}(U), f, g)
\subset (q_1, \dots , q_{n-3}, f, g)
\end{equation}
so by (\ref{eq2linksings}) we have that 
$I_{D_T}(U)\subset I_{{Z_2}}(U)$ and we are done.
Let us consider (\ref{eq3linksings}); since
$\I_{B}(U)\subset
\I_{Y_B}(U)$,
then by (\ref{eq3linksings}) we have
\begin{equation}
\label{eq5linksings}
\I_{D_T}(U)\cdot  \I_{B}(U)
\subset (q_1, \dots , q_{{n-2}\choose 2}, t).
\end{equation}
By (\ref{eq5linksings}) we can write an element of
$\I_{D_T}(U)\cdot  \I_{B}(U)$ as
$\sum_{j=1}^{n-3}h_jq_j+\sum_i h_iq_{n-3+i}+ht$.
Since $t, q_1, \dots, q_{n-3}\in \I_{B}(U)$, then
$\sum_i h_iq_{n-3+i}\in \I_{B}(U)$.
Moreover $\sum_i h_iq_{n-3+i}\in \I_{W}(U)$.
 This implies
$\sum_i h_iq_{n-3+i}\in \I_{W}(U)\cap \I_{B}(U)$.
Since  $B$ and $W$ are geometrically linked by $A$ (this follows from
Remark \ref{rem3})
then $\I_{W}(U)\cap \I_{B}(U)=
\I_{A}(U)=(q_1, \dots, q_{n-3})$. This implies
\begin{equation}
\label{eq6linksings}
\I_{D_T}(U)\cdot  (\I_{B}(U), f, g)
\subset (q_1, \dots , q_{n-3}, f, g, t).
\end{equation}
Since
$(t)\subset I_{Z_3}(U)=(q_1, \dots, q_{n-3}, f, g)$, (\ref{eq6linksings})
 is exactly (\ref{eq4linksings}).

It is so proved that for $i+n-4<\min\{a, b\}$:
\begin{equation}
\label{eq7linksings}
h^0(\I_{{Z_3}|A}(i+n-4))\leq h^0(\I_{{Z_2}|W}(i, n-4))
\end{equation}
Inequality (\ref{eq1linksings}) follows
from (\ref{eq1linkAs}) and (\ref{eq7linksings}).

Let us consider now the general case $i<\min\{a, b\}$
(and $i+n-4\geq\min\{a, b\}$) and let
$a=\min\{a, b\}$. In this case $i+n-4=a+j$ with $0\leq j \leq n-5$.
From the exact sequence:
\[
0\to \I_{{Z_2\cup B}|A}(a+j)\to
\I_{{Z_3}|A}(a+j)\to \I_{{Z_3}|{Z_2\cup B}}(a+j)\to 0,
\]
since $Z_3=Z^\prime+Z_2$ we have that
 $\I_{{Z_3}|{Z_2\cup B}}=\I_{{Z^\prime}|B}$
and  we find
\begin{equation}
\label{eq8linksings}
h^0(\I_{{Z_3}|A}(a+j))-h^0(\I_{{Z^\prime}|B}(a+j))
\leq h^0(\I_{{Z_2\cup B}|A}(a+j)).
\end{equation}
The left side of inequality (\ref{eq8linksings}) is equal
by (\ref{eq2linkAs}) to
$h^1(\I_{{Z_1}|W}(a+b-2-i))+h^0(\I_{{Z_A}|A}(a+j))-
h^0(\I_{{Z^\prime}|B}(a+j))$, while the right side
(we just proved it) is less or equal than
$h^0(\I_{{Z_2}|W}(i, n-4))$. Therefore we prove
(\ref{eq1linksings}) if we  prove that
\begin{equation}
\label{eq9linksings}
h^0(\I_{{Z_A}|A}(a+j))=
h^0(\I_{{Z^\prime}|B}(a+j)).
\end{equation}
From the resolution of $\O_{Z^\prime}$ in $\Mod(B)$
(resp. of $\O_{Z_A}$ in $\Mod(A)$) we find that
$h^0(\I_{{Z^\prime}|B}(a+j))=
h^0(\O_B(j))+h^0(\O_B(a+j-b))-h^0(\O_B(j-b))$
(resp. $h^0(\I_{{Z_A}|A}(a+j))=
h^0(\O_A(j))+h^0(\O_A(a+j-b))-h^0(\O_A(j-b))$).
Therefore we prove (\ref{eq9linksings}) if we prove
$h^0(\O_A(k))=h^0(\O_B(k))$ for every $k\leq n-5$. Since $A$
and $B$ are arithmetically Cohen-Macaulay this is equivalent
to prove that
$h^0(\I_{A|\PP}(k))=h^0(\I_{B|\PP}(k))$,
that is $h^0(\I_{B|A}(k))=0$, for every $k\leq n-5$.
We prove this with a simple calculus of degrees:
every hypersurface of degree $k\leq n-5$ containing $B$
intersects $W$ at least in $Y_B$, since $\deg(Y_B)
=(n-4)(n-3)>k\cdot\deg(W)=k(n-2)$,
the hypersurface contains also $W$, so it contains $A$.
The equivalent statement in terms of the Hilbert function of $Z_1$ can be
now 
be proved exactly like in Prop. \ref{linkss}.
\end{pf}

Let us consider now a singular rational normal
$3$-fold $X\subset \PP^n$.

\begin{theorem}
\label{linksingt}
Let $X\subset\PP^n$ be a singular rational normal $3$-fold.
Let $Y_1, Y_2\subset X$ be  projective $1$-dimensional schemes. Assume
that $Y_1$ is
arithmetically Cohen-Macaulay. Let $Y=X\cap F_a\cap F_b$ a
c.i. of type $(a, b)$ on $X$.
Assume that $Y_1$ and $Y_2$ are geometrically linked
by $Y$.
 Then for $i<\min\{a, b\}$:
\begin{equation}
\label{eq1linksingt}
h^0(\I_{{Y_2}|X}(i, n-4))\geq h^1(\I_{{Z_1}|W}(a+b-2-i)). 
\end{equation}
\end{theorem}
\begin{pf}
Let $A\subset \PP^n$, $B\subset A$ and $Y_B\subset X$ 
as in Lemma \ref{prel} and
Remark \ref{genA}. 
Let $Y_3$ be the
(geometrically) linked scheme to $Y_1$ by the complete intersection
$Y_A=A\cap F_a\cap F_b$. This implies that $Y_3$  is 
arithmetically Cohen Macaulay. Let us suppose $i+n-4<\min\{a, b\}$.
A hypersurface $T$ of $\PP^{n}$ of degree $i+n-4$
containing $Y_3$ but not $A$ cuts on $X$ a
divisor ${(T)}\sim (i+n-4)H$ which splits in the union of $Y_B$
and a divisor $D_T\sim iH+(n-4)R$. $Y_2$ is of course contained
in ${(T)}$,
we claim that $Y_2$ is contained in $D_T$. But this follows
from the fact that a
 general hyperplane section $D_T^H$ of $D_T$ contains a general
hyperplane
section $Z_2$ of $Y_2$, by Th. \ref{linksings}.
So it is proved that for $i+n-4<\min\{a, b\}$:
\begin{equation}
\label{eq2linksingt}
h^0(\I_{{Y_2}|X}(i, n-4))\geq h^0(\I_{{Y_3}|A}(i+n-4)).
\end{equation}
Since $Y_3$ and $A$ are both arithmetically Cohen Macaulay, then
we have that $h^1(\I_{{Y_3}|A}(k))=0$ for every $k$.
Let $Z_3$, $A_H$ and $B_H$ be
general hyperplane sections of $Y_3$, $A$ and $B$ respectively, then
\begin{equation}
\label{eq3linksingt}
h^0(\I_{{Y_3}|A}(k))\geq h^0(\I_{{Z_3}|A_H}(k))
\end{equation}
for every $k$.
 For $i+n-4<\min\{a, b\}$ inequality (\ref{eq1linksingt}) follows
from (\ref{eq2linksingt}), (\ref{eq3linksingt}) and (\ref{eq1linkAs}).

Let us consider now the general case $i<\min\{a, b\}$
(and $i+n-4\geq\min\{a, b\}$) and let
$a=\min\{a, b\}$. In this case $i+n-4=a+j$ with $0 \leq j\leq n-5$.
Let $Y^\prime=B\cap F_a\cap F_b$. In the same way we have proved
(\ref{eq8linksings}) we find:
\begin{equation}
\label{eq4linksingt}
h^0(\I_{{Y_3}|A}(a+j))\leq
h^0(\I_{{Y^\prime}|B}(a+j))+h^0(\I_{{Y_2\cup B}|A}(a+j)).
\end{equation}
Since both $Y^\prime$ and $B$ are arithmetically
Cohen-Macaulay we have $h^1(\I_{{Y^\prime}|B}(k))=0$ for every $k$,
 then $h^0(\I_{{Y^\prime}|B}(a+j))=
h^0(\I_{{Y^\prime}|B}(a+j-1))+
h^0(\I_{{Z^\prime}|{B_H}}(a+j))$. Moreover,
as we just proved,
$h^0(\I_{{Y_2\cup B}|A}(a+j))\leq h^0(\I_{{Y_2}|X}(i, n-4))$.
Therefore (\ref{eq4linksingt}) becomes
\begin{eqnarray}
h^0(\I_{{Y_3}|A}(a+j))&\leq & h^0(\I_{{Y^\prime}|B}(a+j-1))+
h^0(\I_{{Z^\prime}|{B_H}}(a+j))\nonumber \\
{}&+& h^0(\I_{{Y_2}|X}(i, n-4)).\label{eq5linksingt}
\end{eqnarray}
Since $Y_3$ and $A$ are both arithmetically Cohen-Macaulay
we have $h^0(\I_{{Y_3}|A}(k))=h^0(\I_{{Z_3}|{A_H}}(k))
+h^0(\I_{{Y_3}|A}(k-1))$ for every $k$ and therefore
by (\ref{eq2linkAs}) the left side of
inequality (\ref{eq5linksingt}) is
\begin{equation}
\label{eq6linksingt}
h^1(\I_{{Z_1}|W}(a+b-2-i))+h^0(\I_{{Z_A}|{A_H}}(a+j))
+h^0(\I_{{Y_3}|A}(a+j-1))
\end{equation}
In the same way we have proved 
(\ref{eq9linksings})  we find that
$h^0(\I_{{Y_A}|A}(a+j))=
h^0(\I_{{Y^\prime}|B}(a+j))$ (for $j\leq n-5$).
Since $Y_3\subset Y_A$, then
$h^0(\I_{{Y_3}|A}(a+j-1))\geq h^0(\I_{{Y_A}|A}(a+j-1))$, therefore:
\begin{equation}
\label{eq7linksingt}
h^0(\I_{{Y^\prime}|B}(a+j-1))-h^0(\I_{{Y_3}|A}(a+j-1))\leq 0.
\end{equation}
By (\ref{eq5linksingt}), (\ref{eq6linksingt}), (\ref{eq9linksings})
and (\ref{eq7linksingt}) we get
(\ref{eq1linksingt}).
\end{pf}
\begin{corollary}
\label{corlinksingt}
If  the vertex of $X$ is a point, then for $i<\min\{a, b\}$:
\begin{equation}
\label{linksingtp}
h^0(\I_{{Y_2}|X}(i, n-4))\geq h^0(\I_{{Z_2}|W}(i, n-4).
\end{equation}
\end{corollary}
\begin{pf}
In this case a general hyperplane section $W$ of $X$ is smooth.
Then (\ref
{linksingtp}) follows from
(\ref {eq1linksingt}) and (\ref{eq1linkss}).
\end{pf}

At this point we are able to prove the same genus formula
we  have proven in Prop. \ref{genuss} in the case $X$ is singular
with vertex  a  point.

\begin{theorem}
\label{genussing}
Let $Y_1$ and $Y_2$ be  two $1$-dimensional projective schemes
contained in a  rational normal $3$-fold $X\subset\PP^n$ which is
singular whose vertex is  a point.
Let $Y_1$ be locally Cohen-Macaulay.
Assume that $Y_1$ and $Y_2$ are algebraically linked by a complete
intersection
$Y=X\cap F_a\cap F_b$ of type $(a, b)$ on $X$.
Then:
\begin{equation}
\label{genusform2}
p_a(Y_2)=p_a(Y_1)-p_a(Y)+(a+b-3)\cdot\deg(Y_2)+
(n-4)\cdot\deg(R\cap Y_2)+1.
\end{equation}
\end{theorem}
\begin{pf}
Let $A=Q_1\cap \cdots \cap Q_{n-3}$, $B$ and $Y_B=X\cap B$
be as in Lemma \ref{prel}.
Let $Y_A=A\cap F_a\cap F_b$ and let $Y_3$ be  the linked scheme
to $Y_1$
by $Y_A$; let $Y^\prime=Y_3-Y_2$ (as Weil divisor on $A\cap F_a$)
be the linked scheme to $Y_2$ by $Y_3$.
 By adjunction formula $\o_{Y_A}\cong\O_{Y_A}(a+b+n-7)$.
Then with the same techniques used in the proof of Prop. \ref{genuss} 
we obtain:
\begin{equation}
\label{g1}
p_a(Y_3)=p_a(Y_1)-p_a(Y_A)+(a+b+n-7)\cdot \deg(Y_3)+1
\end{equation}
and
\begin{equation}
\label{g2}
p_a(Y^\prime)=p_a(Y)-p_a(Y_A)+(a+b+n-7)\cdot \deg(Y^\prime)+1.
\end{equation}

We claim now that $Y_2$ intersects $Y^\prime$ in the $0$-dimensional
scheme $Y_2\cap Y_B$. Let $U$ be an open subset of $\PP^n$, let $q_i$
for $i=1, \dots, {{n-2}\choose {2}}$
be the equations in $U$ of the quadrics in the homogeneous ideal
$I_{X|\PP}$ of $X$
and let $f_a, f_b$ be the local equations of the hypersurfaces $F_a$ and $F_b$
respectively.
Let
$\I_{{Y_2}}(U)=(q_1, \dots , q_{{n-2}\choose 2}, f_1, \dots, f_l)$ be
the  ideal in $U$ of ${Y_2}_{|U}$, where $f_1, \dots, f_l$
cut $Y_2$ on $X\cap U$. Since $Y_2\subset Y$, then
\[
(q_1, \dots , q_{{n-2}\choose 2}, f_a, f_b)\subset
(q_1, \dots , q_{{n-2}\choose 2}, f_1, \dots, f_l).
\]
The ideal of ${Y^\prime}_{|U}$ in $U$ is 
$\I_{Y^\prime}(U)=(\I_B(U), f_a, f_b)$, therefore we obtain
\begin{align*}
\I_{{Y_2}\cap {Y^\prime}}(U) & =
(q_1, \dots , q_{{n-2}\choose 2}, f_1, \dots, f_l,
\I_{B}(U), f_a, f_b)\\
& = (q_1, \dots , q_{{n-2}\choose 2}, f_1, \dots, f_l,
\I_{B}(U))= \I_{{Y_2}\cap {Y_B}}(U),
\end{align*}
since $\I_{{Y_B}}(U)=
(q_1, \dots , q_{{n-2}\choose 2}, \I_{B}(U))$.

By Noether's formula we  find:
\begin{equation}
\label{g3}
p_a(Y_3)=p_a(Y_2)+p_a(Y^\prime)+\deg(Y_2\cap Y_B)-1.
\end{equation}
Eliminating $p_a(Y^\prime)$ and $p_a(Y_3)$ from (\ref{g3})  using
(\ref{g1}) and (\ref{g2})   we find
\begin{equation}
\label{g4}
p_a(Y_2)=p_a(Y_1)-p_a(Y)+(a+b+n-7)\deg(Y_2)-\deg(Y_2\cap Y_B)+1.
\end{equation}
Since $Y_B\sim (n-4)H-(n-4)R$ by Lemma \ref{prel},
(\ref{g4}) gives exactly (\ref{genusform2}).
\end{pf}

\section{An application to the classification of curves of
maximal genus}

In this section we will show  some examples of application
of the techniques developed in the previous sections to the classification
of curves of maximal genus $G(d, n,  s)$ in $\PP^n$.
Let us first  summarize some results of
\cite{ccd} and some other preliminary  facts
useful to introduce
the problem.

From now on, let $C$ be an integral, nondegenerate curve of degree $d$ and
arithmetic genus $p_a(C)$ in $\PP^n$, with 
$d>\frac{2s}{n-2}\Pi_{i=1}^{n-2} {((n-1)!)}^
{\frac{1}{n-1-i}}$ and $s\geq n-1$ (later we will assume $s\geq 2n-1$).
Assume
$C$ not
contained on surfaces of degree $<s$ and define $m, \epsilon, w, v, k,
\delta$ as follows:

divide $d-1=sm+\epsilon$, $0\le \epsilon\le s-1$ and 
$s-1=(n-2)w+v$, $v=0, \dots , n-3$;

if $\epsilon< w(n-1-v)$, divide $\epsilon=kw+\delta$, $0\le\delta<w$;

if $\epsilon \ge w(n-1-v)$, divide $\epsilon+n-2-v=k(w+1)+\delta$, $0\le
\delta<w+1$.

It is a result of \cite{ccd} (section 5) that the genus $p_a(C)$ is
bounded by the function:
\[
G(d, n, s)=1+\frac{d}{2}(m+w-2)-\frac{m+1}{2}(w-3)+\frac{vm}{2}(w+1)+\rho
\]
where $\rho=\frac{-\delta}{2}(w-\delta)$ if $\epsilon<w(n-1-v)$ and
$\rho=\frac{\epsilon}{2}-\frac{w}{2}(n-2-v)-\frac{\delta}{2}(w-\delta+1)$
if $\epsilon\ge w(n-1-v)$.

If $Z$ is a general hyperplane section of $C$ and $h_Z$ is the Hilbert
function of $Z$, then the difference $\Delta h_Z$ must be bigger than the
function $\Delta h$ defined by:
\[
\Delta h(r)=
\begin{cases}
0 & \text{if} \quad r<0 \cr
(n-2)r+1 & \text{if} \quad  0\le r \le w \cr
s & \text{if} \quad  w< n \le m \cr
s+k-(n-2)(r-m) & \text{if} \quad m< r \le m+\delta \cr
s+k-(n-2)(r-m)-1 & \text{if} \quad m+\delta < r \le m+w+e \cr
0 & \text{if} \quad r>m+w+e 
\end{cases}
\]
where $e=0$ if $\epsilon<w(n-1-v)$ and $e=1$ otherwise (\cite{ccd} Prop.
0.1).
\begin{proposition}
\label{extcurve}
If $p_a(C)=G(d, n, s)$, then $C$ is arithmetically Cohen-Macaulay and
 $\Delta h_Z(r)=\Delta h(r)$ for all $r$.
Moreover $Z$ is contained on a
reduced curve $\Gamma$ of degree $s$ and maximal genus 
$G(s,
n-1)={w\choose 2}+wv$ in $\PP^{n-1}$ 
(Castelnuovo curve). Since
$d>s^2$, $\Gamma$ is unique and, when we move the hyperplane, all these
curves $\Gamma$'s patch togheter giving a surface $S\subset \PP^n$ of
degree
$s$ through $C$ (Castelnuovo surface).
\end{proposition}
\begin{pf} 
See \cite{ccd} Prop. 6.1, Prop. 6.2 and Cor. 6.3.
\end{pf}

\begin{proposition}
\label{surf}
The surface $S$ of Prop. \ref{extcurve}
is irreducible and if $s\geq 2n-1$ it  lies on a 
 rational normal  $3$-fold $X\subset\PP^n$. As a divisor
on $X$ the surface $S$ is linearly equivalent to $(w+1)H-(n-3-v)R$
(or $wH+R$ if $v=0$).
If $n=6$ and $s$ is even  there is the further possibility that the
surface $S$ lies in
a
cone
over the Veronese surface in $\PP^5$ and is the complete intersection  
with a hypersuface not containing the
vertex.
 \end{proposition}
\begin{pf}
$S$ is irreducible since $C$ is irreducible and is not contained on
surfaces of degree $<s$. 
The rest of the statement follows using the characterization of
Castelnuovo surfaces given in \cite{ha2}.

\end{pf}

\begin{proposition}
\label{m+1}
There exists a hypersurface
$F_{m+1}$ of degree $m+1$, passing through $C$ and not containing $S$.
\end{proposition}
\begin{pf}
For a general hyperplane section $\Gamma$ of $S$, the Hilbert function
$h_\Gamma$ is known (see e.g. \cite{ha1} Th. 3.7); in particular we have 
$\Delta
h_\Gamma (r)=\Delta h_Z (r)$ when $0\le r\le m$ and hence
$h^0({\I}_C (r))=h^0({\I}_S(r))$ when $0\le r\le m$. For $r=m+1$
one computes $\Delta h_\Gamma (m+1)<\Delta h_Z (m+1)$ and this implies
$h^0({\I}_{C|{\PP}} (m+1))>h^0({\I}_{S|{\PP}}(m+1))$.
\end{pf}

Let us suppose $s\geq 2n-1$. By the  Prop. \ref{surf}
a  curve $C\subset\PP^n$ of maximal genus $G(d, n, s)$ 
lies then on a rational normal
$3$-fold $X$ (except in the case where $S$ lies in a cone over a
Veronese surface, which we do  not intend to go through). Let $F_{w+1}$
be a hypersurface of degree $w+1$ cutting $S$ on $X$. By Prop. \ref{m+1} 
we can consider  on $S$ the curve $C^\prime$ 
 residual to $C$ by the intersection
with the hypersurface $F_{m+1}$. Since $\deg(C^\prime)<\deg(C)$, then
$C^\prime$ does not contain $C$.  
Choosing in $X$ a sufficiently general  divisor  $D\sim (n-3-v)R$
(or $D\sim H-R$ in case $S\sim wH+R$) linked  
to $S$ by  $X\cap F_{w+1}$, then the
residual  scheme
on $X$ to $C$  by
the c.i.
$X\cap F_{w+1}\cap F_{m+1}$ is a curve which we call
$C^{\prime\prime}$. 
When $v=n-3$ then $S=X\cap F_{w+1}$
and of course $C^\prime=C^{\prime\prime}$; otherwise $C^{\prime\prime}$
is the union of $C^\prime$ with a curve $C_D$ 
contained in $D$, therefore
$C_D$ is formed by
$n-3-v$ distinct plane
curves 
of degree $m+1$ or, in case $S\sim wH+R$, $C_D$ is the complete
intersection on  $D\sim H-R$ by a hypersurface of degree $m+1$. Letting
$Z^\prime$, $Z^{\prime\prime}\subset W$ be
general
hyperplane
sections of $C^\prime$ and $C^{\prime\prime}$ respectively, 
we have the following two Lemmas:

\begin{lemma}
\label{first}
If $X$ is smooth or if the vertex of $X$ is a point, then for $i\leq w, m$
\[
h^0(\I_{C^{\prime\prime}|X}(i, n-4))\geq
h^0(\I_{Z^{\prime\prime}|W}(i, n-4)) =
\sum_{r=m+w-i+1}^\infty \Delta h(r).
\]
Moreover if $h^0(\I_{Z^{\prime\prime}|W}(i-1, n-4))=0$ and
$h^0(\I_{Z^{\prime\prime}|W}(i, n-4))=h>0$, then 
$h^0(\I_{C^{\prime\prime}|X}(i-1, n-4))=0$
and $h^0(\I_{C^{\prime\prime}|X}(i, n-4))=h$.
\end{lemma}
\begin{pf}
$C$ and $C^{\prime\prime}$ are geometrically linked
by $Y=X\cap F_{w+1}\cap F_{m+1}$ since they are equidimensional, have no 
common components ($C$ is irreducible and $C^\prime$ does not contain $C$)
and  no embedded components
($Y$ is arithmetically Cohen Macaulay). $W$ is smooth and $C$
is arithmetically Cohen Macaulay, therefore
by Prop. 
\ref{linkss} we know that 
$
h^0(\I_{Z^{\prime\prime}|W}(i, n-4))= 
d-h_Z(m+w-i)
$. Then  note that for every $k$ 
we have $d-h_Z(k)=d+\Delta h_Z(k+1)- h_Z(k+1)=d+\sum_{r=k+1}^t \Delta
h_Z(r)-h_Z(t)
=\sum_{r=k+1}^\infty \Delta h_Z(r)$ because for $t$ big we
have $h_Z(t)=d$, and that, by Prop. \ref{extcurve},
$\Delta h_Z(r)=\Delta h (r)$ for all $r$.
By Cor. \ref{acm} (if $X$ is smooth) or by
Cor. \ref{corlinksingt} (if the vertex of $X$ is a point) we have that 
\begin{equation}
\label{bigger}
h^0(\I_{C^{\prime\prime}|X}(i, n-4))\geq
h^0(\I_{Z^{\prime\prime}|W}(i, n-4)). 
\end{equation}
From the exact sequence
\[
0\to \I_{C^{\prime\prime}|X}(k-1, n-4)\to
\I_{C^{\prime\prime}|X}(k, n-4)\to
\I_{Z^{\prime\prime}|W}(k, n-4)\to 0
\]
we obtain that if for $k=i-1$ we have $h^0(\I_{Z^{\prime\prime}|W}(k,
n-4))=0$, then $h^0(\I_{C^{\prime\prime}|X}(k, n-4))=0$. In this
hypothesis for $k=i$ we have an injection
$H^0(\I_{C^{\prime\prime}|X}(i, n-4)) \hookrightarrow 
H^0(\I_{Z^{\prime\prime}|W}(i, n-4))$ and therefore 
by 
(\ref{bigger})
 $h^0(\I_{C^{\prime\prime}|X}(i,
n-4)) = h^0(\I_{Z^{\prime\prime}|W}(i, n-4))$.
\end{pf}
\begin{lemma}
\label{second}
If the vertex of $X$ is a line, then for $i\leq w, m$
\[
h^0(\I_{C^{\prime\prime}|X}(i, n-4))\geq
h^1(\I_{Z|W}(m+w-i))=\sum_{r=m+w-i+1}^\infty h(r).
\]
\end{lemma}
\begin{pf}
Th. \ref{linksingt}.
\end{pf}

The strategy is to classify all the curves of maximal genus in $\PP^n$ for
arbitrary $n$ by classifing
the  linked curves $C^\prime$'s. A complete classification
Theorem when $n=4$ is proved in \cite{cc} and  when $n=5$
in \cite{f1} (and in the forthcoming work \cite{f3}). Depending on the numerical
parameters ( $\epsilon$, $w$, $v$, $k$) associated to $C$  and
on the
type of the scroll $X$ the analysis goes on case by case.
In the following example we want to show the simplest non trivial case in
the classification
procedure, when $C^\prime$ is a plane curve (the trivial case is
$C^\prime=\emptyset$).   It should be remarked that while in $\PP^3$ the curve
$C^\prime$ is always degenerate this is no longer true for $n\geq 4$
(see \cite{cc} and \cite{f1} for $n=4, 5$).

\begin{example}
\label{esempio1}
Let $s\geq 2n-1$ and let $d>\frac{2s}{n-2}\Pi_{i=1}^{n-2} {((n-1)!)}^
{\frac{1}{n-1-i}}$; divide $s-1=(n-2)w+v$, $v=0, \dots
, n-3$ and divide $d-1=sm+\epsilon$, $0\leq\epsilon\leq s-1$. 
Suppose $s-2-w\leq\epsilon\leq s-2$. Let $C\subset\PP^n$ be a curve 
of maximal genus $G(d, n, s)$. Then  the linked curve $C^\prime$
 is a plane curve
 of degree $s-\epsilon-1$.
In case that the vertex of $X$ is  a line, $C^\prime$ will not contain
this line
as a component.
\end{example}

Here  we suppose for the sake of  simplicity  that $v=n-3$ (the
result can be proved with similar arguments
for every $v$), i.e. we
put
ourselves in the simplest case $C^\prime=C^{\prime\prime}$; with this
assumption we always have $e=1$, i.e.
$\epsilon\geq w(n-1-v)=(n-3)(w+1)$, hence we write
  $\epsilon+1=k(w+1)+\delta$ with 
$k=n-3$ and $\delta\leq w$.  
In this case $C$
and $C^\prime$ are (geometrically)
linked  by a c.i. $Y=X\cap F_{w+1}\cap F_{m+1}$
on $X$.
If $X$ is smooth or if the vertex of $X$ is a point, then applying Lemma
\ref{first}
for $i=0$ (of course $h^0(\I_{C^{\prime}|X}(-1, n-4))=0$)
we compute:
\[
h^0(\I_{C^{\prime}|X}(0, n-4))=
n-4.
\]
The linear system $|\O_X(0, n-4)|$ is composed with a rational 
pencil,  i.e. we have $\pi: X\to \PP^1$ and 
$|\I_{C^\prime|X}(0, n-4))|=\pi^* \G$, where $\G$ is a linear subsystem
of
$|\O_{\PP^1}(n-4)|$.
This implies that 
$|\I_{C^{\prime}|X}(0, n-4)|$ has a fixed part;
in this case, since $h^0(\O_X(0, a))=a+1$ for every $a\geq 0$, the
fixed part  of $|\I_{C^\prime}(0, n-4)|$ is $\sim R$ and the moving part
is equal to the whole  
$|\O_X(0, n-5)|$. 
Therefore we  conclude that
$C^\prime$ is contained in a plane $\pi\sim
R$.


If the vertex of $X$ is a line 
 then applying Lemma \ref{second}  we compute:
\[
h^0(\I_{C^{\prime}|X}(0, n-4))\geq
n-4.
\]
Let us suppose for the moment that
 $h^0(\I_{C^{\prime}|X}(0,
n-4))=n-4$ (we will  exclude the case
$h^0(\I_{C^{\prime}|X}(0, n-4))>n-4$ in the sequel).
We want to   conclude as in the previous case  
 that
$|\I_{{C^\prime|X}}(0, n-4)|$ has a fixed part. So let
us suppose that  
$|\I_{C^{\prime}|X}(0, n-4)|$ has no fixed part, which implies
that the support of $C^\prime$ is the singular line of $X$.
By Bertini's Theorem  the
generic divisor in the corresponding linear subsystem $\G$ of $\PP^1$ is 
union of $n-4$ distinct points in a rational normal  curve $C_{n-2}$ of
degree $n-2$, which span a $\PP^{n-5}$. 
Therefore we can choose a basis $\{D_1,
\dots 
, D_{n-4}\}$ in the linear system $|\I_{C^{\prime}|X}(0, n-4)|$ such
that $D_i$ is union
of $n-4$ distinct planes of $X$ for every $i$ and
such that the linear space spanned by  each $D_i$ is
$<D_i>\cong\PP^{n-3}$ and $D_i=X\cap <D_i>$. In this situation 
the base locus
of $|\I_{C^{\prime}|X}(0, n-4)|$,
which is equal to $D_1\cap \cdots \cap D_{n-4}=X\cap <D_1>\cap \cdots \cap
<D_{n-4}>$, is necessarily the singular line $l\cong \PP^1$ of
$X$
counted
with multiplicity one, but this is not possible since
$|\I_{l|X}(0, n-4))|= |\O_X(0, n-4)|$ and we  have a contradiction.
Therefore, as in the previous case, we conclude that $C^\prime$ is contained in a plane $\pi\sim R$.
We claim now that $C^{\prime}$ cannot contain the singular line
of $X$ as a component. In fact in this case both $S$ and $F_{m+1}$ would
pass through it and their proper transforms $\tilde S$ and $\tilde
{F}_{m+1}$ on the canonical resolution $\tilde X$ of $X$ would be
$\tilde S\sim (w+1-a)\tilde H+ (n-2)a\tilde R$ and $\tilde {F}_{m+1}\sim
(m+1-b)\tilde H+(n-2)b\tilde R$ with $a, b\geq 1$.
In this case, since $C$ is irreducible (therefore it does not contain the
singular line), $C^\prime$ would contain the singular line with
multiplicity $\alpha$ which we compute using \cite{f2} Prop. 4.11 as:
$\alpha=S^*\cdot F^*_{m+1}\cdot\tilde H-\tilde S\cdot \tilde
{F}_{m+1}\cdot
\tilde H
=(w+1)\tilde
H\cdot(m+1){\tilde
H}^2-\tilde S\cdot \tilde {F}_{m+1}\cdot
\tilde H=
ab(n-2)$, where $S^*$ and $F^*_{m+1}$ are respectively the integral total
transform of $S$ and $F_{m+1}$ (Def. \ref{totale}).  
But since  $C^\prime$ is contained in a plane
$\pi\sim 
R$ by the same kind of computation we conclude that 
$C^\prime$ would contain the singular line
with multiplicity
$\beta=S^*\cdot R^*\cdot\tilde H-\tilde S\cdot \tilde R\cdot
\tilde H =
(w+1)\tilde
H\cdot (\tilde H
-(n-3)\tilde R)\cdot\tilde
H-((w+1-a)\tilde H+(n-2)a\tilde
R)\cdot \tilde R\cdot \tilde H=a$ and this is in contradiction with the
previous value.

In case   
$h^0(\I_{C^{\prime}|X}(0, n-4))>n-4$ we would 
have $h^0(\I_{C^{\prime}|X}(0, n-4))=n-3$, i.e.
$|\I_{C^{\prime}|X}(0, n-4))|=|\O_X(0, n-4)|$.
This would
 imply that $C^\prime$ has
degree one and coincides with the singular line of $X$ but this can be
excluded
with the above computation.
\bigskip
In the next example we show that in the  case of Example \ref{esempio1} 
{\it smooth} curves of maximal genus do always exist. 
Moreover we
explicity construct
such
curves  on a smooth rational normal $3$-fold where 
 genus formula (\ref{genuss}) holds. It is interesting to note that it is
not always possible to construct curves of maximal genus on a {\it smooth}
rational normal $3$-fold. There are cases (for some values of
$d$
and $s$)
where the construction is possible only on  a rational normal $3$-fold
whose vertex is a point and where 
the genus formula (\ref{genussing}) holds, as showed in
 \cite{f1} (Prop. 4.2 part 4) and Example 5.2 case $k=v=1$)
 for $n=5$. The existence of
curves of maximal genus in $\PP^5$ is proved for all cases in \cite{f1}.
We state first the following, easy to prove, result (see \cite{ro} Lemma 1 pg. 133) which we will use later.

\begin{lemma}
\label{rogora}
Let $X$ be a smooth $3$-fold. Let $\Sigma$ be a linear system of surfaces
of $X$ and let $\gamma$ be a curve contained in the base locus of $\Sigma$.
Suppose that the generic surface of $\Sigma$ is smooth at the generic 
point of $\gamma$ and that it has at least a singular point which is variable 
in $\gamma$. Then all the surfaces of $\Sigma$ are tangent along $\gamma$.
\end{lemma}

\begin{example}
\label{esempio2}
For every $d$ and $s$ in the range 
 of $Example$ \ref{esempio1} there exists a smooth curve
$C\subset \PP^n$ of maximal genus $G(d, n, s)$.
\end{example}
 For the sake of simplicity we treat only the cases $v=n-3, n-4$, i.e.
$s=(n-2)(w+1)$ 
and
$s=(n-2)w+n-3$. The other cases can be treated in a similar way.

Let us suppose  $v=n-3$.
Let $X\subset \PP^n$ be a smooth rational normal $3$-fold of degree
$n-2$ and
let $\pi\sim R$ be a plane contained in $X$.
Let $D$ be a smooth curve on $\pi$ of degree \,
$0\leq \deg D=w+1-s+\epsilon+1=\epsilon+1-(n-3)(w+1)\leq w$
(possibly $D=\emptyset$). If we consider the union of $D$
with any plane curve $C^\prime\subset\pi$ of degree 
$w+1-\deg D=s-\epsilon-1$, then
 there exists a hypersurface $F_{w+1}$ 
of degree $w+1$ cutting $C^\prime\cup D$ on $\pi$.
Therefore the linear system
$|\I_{{D}|X}(w+1)|$ 
of divisors on $X$ cut by hypersurfaces of
degree $w+1$ through
$D$ is not empty and cut on $\pi$ the linear system 
$D+|\O_\pi(s-\epsilon-1)|$.  Moreover the linear system 
$|\I_{{D}|X}(w+1)|$ contains  
 the linear
subsystem $L+|\O_X(w)|$, where $L$ 
is a fixed hyperplane section
containing $\pi$, that has fixed part $L$ and no other base points.
This implies that
$D$
is  the base locus of all $|\I_{{D}|X}(w+1)|$ and that $|\I_{{D}|X}(w+1)|$
is not composed with a pencil, because in this case every element in the 
system would be a sum of algebraically equivalent divisors, while the
divisors in  $L+|\O_X(w)|$ are obviously not of this type.
 By Bertini's
Theorem we can then  conclude that the generic divisor in
$|\I_{{D}|X}(w+1)|$ is an irreducible 
surface $S$
of degree
$(n-2)(w+1)=s$   smooth outside
$D$. We claim that  $S$ is in fact smooth at every point $p$ of $D$.
To see this, by Lemma \ref{rogora}, it is enough to prove that, for every  
$p\in D$, there exists a surface in $|\I_{D|X}(w+1)|$ which is smooth at
$p$, and that for a generic point $q\in D$, there exist two surfaces in
$|\I_{D|X}(w+1)|$ with distinct tangent planes at $q$.
In fact, for every $p\in D$ we can always find a surface $T$ in the
linear system $|\O_X(w)|$ which does not pass through $p$, therefore
the surface $L+T$ is smooth at $p$ with tangent plane $\pi$. Moreover a
generic surface in the linear system $|\I_{D|X}|$ which cut $D$ on
$\pi$ has at $p$ tangent plane $T_p\neq \pi$.

Let $\ci\subset\pi$ be the linked curve to $D$ by the intersection
$\pi\cap S$. 
Let us consider the linear system $|\I_{{C^\prime}|S}(m+1)|$ of
divisors cut on $S$ by the hypersurfaces of degree $m+1$ passing
through
$C^\prime$. With the same argument used above we conclude that
 this linear system is not composed with a pencil, it has $C^\prime$
as a  fixed part  
 and no other base points. 
Therefore by Bertini's theorem we deduce that the generic  curve
$C=S\cap F_{m+1}-C^\prime$ in the movable part of 
the linear system is irreducible, smooth 
and has
the required degree
$d=s(m+1)-s+\epsilon+1$. 
By Clebsch formula one computes:
\[
p_a(C^\prime)=\frac{1}{2}((n-2)w+n-4-\epsilon)((n-2)w+n-5-\epsilon).
\]
Moreover $\deg(R\cap C^\prime)=0$. Substituting these expressions in the 
genus
formula (\ref{genuss}) we find that $p_a(C)$ has the maximal value $G(d,
n, s)$, therefore $C$ is the
required curve.

We consider now the case $v=n-4$.
 Let $\pi\sim R$  and  $p\sim R$ be two distinct planes contained in $X$.
Let $D$  be a smooth curve on $\pi$ of degree
\, $0\leq \deg D= \epsilon+2-(n-3)(w+1)\leq w$
(possibly $D=\emptyset$).
Let us consider the linear system
$|\I_{{D\cup p}|X}(w+1)|$ 
 of divisors on $X$ cut by
hypersurfaces of degree $w+1$ containing the plane $p$ and passing
through $D$. This linear system is not
empty since hypersurfaces which are union of a hyperplane containing the
plane $p$ and of a hypersurface of degree
$w$
passing through $D$ cut
on $X$ divisors in the system.
From this description one can see that
$|\I_{{D\cup p}|X}(w+1)|$ 
is not composed with a pencil and that its 
 base locus is  $p\cup D$. By Bertini's Theorem the generic element
in the movable part
of the linear system is an irreducible surface $S\sim (w+1)H-R$ of degree
$s$, smooth outside $D$. 
By the same  argument used in
the previous case we can prove that $S\cup p$ is smooth at every point of 
$D$, but since $D\cap p =\emptyset$ this means that $S$ is 
smooth at $D$.
Let $\ci\subset\pi$ be the linked curve to $D$ by the intersection
$S\cap \pi$. Let us consider the linear system
$|\I_{{C^\prime}|S}(m+1)|$, which is not empty since $\deg\ci <m+1$ and has
base locus equal to the curve $\ci$.
As in the previous case we deduce that the generic  curve
$C=S\cap F_{m+1}-C^\prime$ in the movable part of 
this linear system is irreducible, smooth 
and has
the required degree $d=s(m+1)-\deg(C^\prime)$.
By generality the hypersurface $F_{m+1}$ does not contain the plane $p$
and cut on it a curve $C_1$ of degree $m+1$.
Let
$C^{\prime\prime}=C^\prime\cup C_1$; by construction the curve
$C^{\prime\prime}$
is
 linked to $C$ by a c.i. on $X$ of type $(w+1, m+1)$
 
By Noether's formula one computes:
\[
p_a(C^{\prime\prime})=\frac{1}{2}((n-2)w+n-5-\epsilon)((n-2)w+n-6-\epsilon)
+\frac{1}{2} m(m-1)-1.
\]
Moreover $\deg(R\cap C^{\prime\prime})=0$. Substituting these expressions
in
the 
genus
formula (\ref{genuss}) we find that $p_a(C)$ has the maximal value $G(d,
n, s)$. Therefore $C$ is the
required curve.

\end{document}